\documentclass[11pt,amsart]{article}

\usepackage{epsf}

\usepackage{amssymb,amsthm,verbatim}
\newtheorem{prop}{Proposition}[section]
\newtheorem{thm}[prop]{Theorem}
\newtheorem{lem}[prop]{Lemma}
\newtheorem{cor}[prop]{Corollary}
\newtheorem{conj}[prop]{Conjecture}

\theoremstyle{definition}

\newtheorem{defn}[prop]{Definition}

\newtheorem{ques}[prop]{Question}

\newcommand{\ep}{\epsilon}

\newcommand{\Ga}{\Gamma}

\newcommand{\cC}{\mathcal C}

\newcommand{\bmat}{\left ( \begin{matrix} }
\newcommand{\emat}{\end{matrix} \right ) }

\newcommand{\ben}{\begin{enumerate}}

\newcommand{\een}{\end{enumerate}}

\newcommand{\sumlim}{\sum\limits}

\newcommand{\blem}{\begin{lem}}

\newcommand{\elem}{\end{lem}}

\newcommand{\bcl}{\begin{clm}}

\newcommand{\ecl}{\end{clm}}

\newcommand{\bthm}{\begin{thm}}

\newcommand{\ethm}{\end{thm}}

\newcommand{\bq}{\begin{ques}}
\newcommand{\eq}{\end{ques}}

\newcommand{\bpr}{\begin{prop}}

\newcommand{\epr}{\end{prop}}

\newcommand{\bco}{\begin{cor}}

\newcommand{\eco}{\end{cor}}

\newcommand{\bcon}{\begin{conj}}

\newcommand{\econ}{\end{conj}}

\newcommand{\bde}{\begin{defn}}

\newcommand{\ede}{\end{defn}}

\newcommand{\bex}{\begin{exa}}

\newcommand{\eexa}{\end{exa}}

\newcommand{\bexe}{\begin{exe}}

\begin{document}

\title{Kazhdan constants of group extensions}

\author{Uzy Hadad}
\date{\today}                   % The date
\maketitle                      % This command creates the title for the document

%======================================================================================================================
\begin{center}
%Draft - Please do not distribute.

\end{center}

\begin{abstract} We give bounds on Kazhdan constants of abelian extensions of
(finite) groups. As a corollary, we improved known results of Kazhdan
constants for some meta-abelian groups and for the relatively free group
in the variety of $p$-groups of lower $p$-series of class $2$.
Furthermore, we calculate Kazhdan constants of the tame automorphism
groups of the free nilpotent groups.

\end{abstract}

\section{Introduction.}

Property (T) was introduced by Kazhdan \cite{Kaz} in 1967. It founds numerous
applications in various areas of mathematics, in particular for
estimating quantitative behavior of the mixing time of random walks and
the expansion of the Cayley graphs of finite groups. For introduction to
the subject and applications we refer the reader to \cite{BHV,Lub}.

\begin{defn} \label{def:almost.invariant} Let $\Ga$ be a discrete group, $S\subset\Ga$ a subset, $\ep>0$,
 and let $(\rho,\mathcal{H}_\rho)$ be
 a unitary representation of $\Ga$. A vector $0 \neq v\in \mathcal{H}_\rho$ is
 called $(S,\ep)$-invariant, if for every $s\in S$, we
  have $\|\rho(s)v-v\|\leq \epsilon \|v\|$.
\end{defn}

\begin{defn} \label{defn:Kazhdan.constant} Let $\Ga$ be a discrete group and let $S$
be a finite generating subset.
\begin{enumerate}
\item The group $\Gamma$ is said to have \textit{Kazhdan property
    (T)} if there exists an $\ep>0$ such that every unitary
    representation with $(S,\ep)$-invariant vector, contains a
    non-zero $\Gamma$-invariant vector. In that case $\ep$ is called
    a \textit{Kazhdan constant} of $\Ga$ with respect to $S$.

\item More generally, if $\cC$ is a class of unitary representations
    of $\Ga$, we say that $\ep$ is a \textit{Kazhdan constant} of
    $\Ga$, with respect to the set $S$ and relative to the class
    $\cC$, if every representation in $\cC$ with $(S,\ep)$-invariant
  vector, contains a non-zero $\Gamma$-invariant vector.

\end{enumerate}
\end{defn}

It is a standard observation that a group $\Gamma$ has T iff
there exists a finite generating multiset $S$ (i.e. an object like a set except that multiplicity is significant) of $\Ga$ and an $\ep>0$ such that

    $$\inf \limits_{\rho \in \Gamma_0^*}
    \inf \limits_{0 \neq v \in \mathcal{H}_{\rho}} \frac{1}{|S|}\sum_{s\in
    S}\frac{\|\rho(s)v-v\|^2}{\|v\|^2} \geq \ep,$$

\noindent where  $\Gamma_0^*$ denote the space of unitary representations
     of $\Gamma$ with no invariant non-zero vectors.

     \bde  The number $\ep$ above is called a {\it average Kazhdan constant} of $\Gamma$ with respect to $S$.
\ede
\begin{comment}
More generally, if $\cC$ is a class of unitary
representations
    of $\Ga$ with no invariant non-zero vectors, then $\ep>0$ is
    called an \textit{average Kazhdan constant} of $\Ga$, with respect to
    the set $S$ and relative to the class $\cC$, if
     $$ \inf \limits_{\rho \in \cC} \inf \limits_{0
    \neq v \in \mathcal{H}_{\rho}} \frac{1}{|S|}\sum_{s\in
    S}\frac{\|\rho(s)v-v\|^2}{\|v\|^2} \geq \ep.$$

\end{comment}

\noindent This definition is equivalent to the previous definition, they only
differ in $\ep$.
 One can check that the infimum in the average Kazhdan constant
is obtained at an irreducible representation.
In general, for a group $\Gamma$ with a Kazhdan constant $\ep_1$ with respect to $S$,
and with an average Kazhdan constant $\ep_2$ (with respect to $S$), we have the
following inequalities $\ep^2_1 \geq \ep_2 \geq \frac{\ep^2_1}{|S|}$.

It turns out (see subsection \ref{spectral_gap}) that in case of a finite group, it is more
appropriate to define a Kazhdan constant of a group with respect
to the mean instead of the maximum over the generators. So throughout this paper the notation \textit{Kazhdan constant}
 will be used for the case of
infinite groups and an  \textit{average Kazhdan constant} for finite groups.

Let $\Gamma$ be a finitely generated group with a Kazhdan constant $\ep_1$ with respect to
$S_1$. Let $S_2$ be another generating set then it is not hard to
prove (see \cite{Lub}, rem. 3.2.5) that there exists
$\ep_2=\ep_2(\ep_1,S_1,S_2)$ such that the Kazhdan constant of
$\Gamma$ with respect to $S_2$ is greater than or equal to
$\ep_2$. For the interesting question of the
dependance of Kazhdan constants on the generating sets, see
\cite{ALW,GZ,Kas1}.
 This bring us to the following definition:

\begin{defn} Let $\Ga$ be a finite group. Given $0<\ep<1$, we define $g_\ep(\Ga)$ as the minimal integer $N$,
such that there is a generating set $S$ of $\Ga$, with $|S|\leq N$, and
the average Kazhdan constant of $\Ga$ with respect to $S$ is greater
 than or equal to $\ep$.
\end{defn}

The number $g_{\ep}(\Gamma)$ is important in many theoretical and
practical computational problems, hence we ask the following
natural question.
 \bq
For a finite group $\Gamma$ and $0< {\ep} <1$, what is $g_{\ep}(\Gamma)$?

\eq

\noindent N. Alon and Y. Roichman (Corollary 1 in \cite{AR}, see also
\cite{LR,LS}) showed that $g_{\ep}(\Gamma) \leq
O_{\ep}(\log|\Gamma|)$ for any $0<{\ep}<1$. This bound is sharp for
the case of abelian groups. However, for any finite non-abelian
simple group $\Gamma$ (except Suzuki), it has been proved recently
by Kassabov, Lubotzky and Nikolov \cite{KLN}, that
$g_{\ep}(\Gamma)=O(1)$ for some $0< \ep <1$.

One of the goals of this paper is to present a method for estimating $g_{\ep}(\Gamma)$
for some non-simple groups.

\subsection{Main Results}

 \bthm \label{thm_gen_ext}

Let $A, H$ be two finite groups with $A$ abelian, such that $H$ acts on
$A$. Let $\Gamma$ be an extension of $A$ by $H$. Let $S \subset H$,
$\tilde{S}$ a lifting of the set $S$ to $\Gamma$, $B \subset A$ a subset and
$\tilde{B}=\{a^x \mid a \in B, x \in H\}$ be the union of the orbits.
Assume that the following hold:

\begin{enumerate}

\item $H$ is generated by  $S$ and has average Kazhdan constants
    ${\ep}_H$ with respect to $S$.

\item $A$ is generated by $\tilde{B}$ and has average Kazhdan
    constants ${\ep}_A$ with respect to $\tilde{B}$.

\end{enumerate}

    Then the group $\Gamma$ has average Kazhdan constants $$ {\ep} > \frac{{\ep}_H
{\ep}_A}{512\cdot\left(1+\frac{|S|}{|B|}+\frac{|B|}{|S|}\right) }$$ with respect to the set of
generators $\tilde{S} \cup B$.

 \ethm

\noindent Note that the conclusion of the theorem holds for an arbitrary lifting
$\tilde{S}$ of $S$, with the same constant. Moreover, the assumptions on
$A$ are with respect to the generating set $\tilde{B}$ which may be very
large (of size $|B|\cdot |H|$) whereas the generating set of $\Gamma$
might be much smaller.

This theorem is an extension of a work by Alon, Lubotzky and Wigderson
\cite{ALW} who proved a similar result for the case of sem-direct product. In fact they showed that
the semi-direct product is a group theoretic analog to the zig-zag product of
Reingold, Vadhan and Wigderson \cite{RVW}.
Here we show that a similar result holds also for other cases of extensions as well.

 The main idea of the proof of Theorem
\ref{thm_gen_ext} is the connection between almost $\Gamma$-invariant
vector and almost $H$-invariant measure on the unitary dual of $A$. If
there exists an almost invariant vector with respect to the generating set of
$\Gamma$, then one can associate with it a measure on the unitary dual of
the abelian group $A$. We show that this measure is almost $H$-invariant
under the dual action of $H$ on $\hat{A}$, and this brings to a
contradiction.

Our method is generic and makes no assumption on the finite groups $A$
and $H$. However, for some specific groups (even infinite) $A$ and $H$,
one can do better. For example, see \cite{Bur,Kas1,Sh1}, where the authors study the
 group $EL_n(R) \ltimes R^n$ (where $R$ is an associative ring) with respect to some unipotent
matrices as a generating set.

The following is deduced from the main Theorem:

\begin{cor}\label{cor_meta_abel}
Let $n \in\mathbb{N}$, $p$ a prime, and let $C_p,C_{p^n}$ be the cyclic
groups of order $p$ and $p^n$ respectively.  Let $\Gamma$ be any
extension of the group algebra $\mathbb{F}_p[C_p]$ by $C_{p^n}$, where $C_{p^n}$ acts by left
translation (mod p) on $\mathbb{F}_p[C_p]$. Then there exists  ${\ep}>0 $, such
that $g_{{\ep}}(\Gamma)\leq O_{{\ep}}(n\log p)$ .

\end{cor}

\noindent Note that by Alon-Roichman Theorem, $g_\ep(\Gamma) \leq O_\ep((n+p)\log
p)$.

\noindent Theorem \ref{thm_gen_ext} covers the case of non-trivial extensions. For
the case of central extension the following is a quantitative version of
a result of Serre.

 \bthm [Serre] \label{thm_cen_ext}

Let $\Gamma$ be a central extension of $A$ by $H$, i.e $A\subseteq
\textbf{Z}(\Gamma)$. Let $\tilde S \subset \Gamma$ be a set of generators
of $\Gamma$ and let $S$ be the projection of $\tilde{S}$ to $H$. Assume
that $\ep_1$ (resp. $\ep_2$) is a Kaz. constant (resp. an average Kaz. constant) of
$H$ with respect to $S$. Then  $\frac{\ep_1}{2\cdot|S|}$ (resp. $\frac{\ep_2}{4}$)
is a Kaz. constant (resp. an average Kaz. constant) of $\Gamma$ with respect to $\tilde{S}$
  with respect to any unitary representation of $\Gamma$ which does not
contains a one dimensional subrepresentation. \ethm

\noindent For the proof we refer the reader to \cite{HV} (pages 26-28 or Lemma
1.7.10 in \cite{BHV}). Although the quantitative result is not stated
there explicitly, it follows immediately from the proof. In subsection
\ref{auto_free_nil} we use this result for estimating Kazhdan
constants of the tame automorphism groups of the free nilpotent groups,
extending a result of Mallahi \cite{MAK}. These estimations have
applications for evaluating the mixing time of the product replacement
algorithm.

To illustrate the use of the above result for concrete groups, we prove
the following:

\noindent For a fixed prime $p$, let $\Gamma_{k,2}$ be  be the relatively free group on $k$ generators
in the variety of $p$-groups of lower $p$-series of class $2$; the precise
definition is provided in Sec. \ref{rel_free}.

\begin{cor} \label{col_free_nil_p} For all $0 < {\ep} < \frac{1}{4}$ we have $$g_{{\ep}}(\Ga_{k,2})\leq
O_{{\ep}}(\sqrt{\log|\Ga_{k,2}|}).$$
\end{cor}
\noindent So our method improves Alon-Roichman result which requires
$O_{\ep}(\log|\Gamma_{k,2}|)$ elements.

\subsection{The Spectral Gap} \label{spectral_gap}
Let $G=Cay(\Gamma,S)$ be the Cayley graph of a finite group
$\Gamma$ with respect to a symmetric generating set $S$
($S=S^{-1}$). Let $A_r$ be the normalized adjacency matrix of the
Cayley graph $G$. We denote by $\lambda(\Gamma,S)$ the second
largest eigenvalue of $A_r$ and we also define its \emph{spectral
gap} $\beta=1-\lambda(\Gamma,S)$ (for applications of $\beta$ see
\cite{HLW}). It turns out that there is a tight connection between
the average Kazhdan constant and the spectral gap of $G$, indeed
we have (see fact 1.1 in \cite{MW} and section 11 in \cite{HLW}):
$$\beta=\inf \limits_{\rho \in \Gamma_0^*} \inf \limits_{0
    \neq v \in \mathcal{H}_{\rho}} \frac{1}{2|S|}\sum_{s\in
    S}\frac{\|\rho(s)v-v\|^2}{\|v\|^2}=
    \inf \limits_{\rho \in \Gamma_0^*} \inf \limits_{0
    \neq v \in \mathcal{H}_{\rho}}\frac{ \|\frac{1}{|S|}\sum_{s\in
    S}{\rho(s)v-v\|}}{\|v\|},$$

\noindent where $\Gamma_0^*$ denotes the space of unitary representations of $\Gamma$ with no invariant non-zero vectors.

 Let us now describe the
structure of this paper. Section \ref{abel_ext} is devoted to a proof
of Theorem \ref{thm_gen_ext}. In subsection \ref{solv_grp} we show that
if $n >1$, then the cohomology group classifying the group extensions of $F_p[C_p]$
by $C_{p^n}$, is not trivial. Moreover, we prove the existence of ``good
orbits" in $F_p[C_p]$ under the action of $C_{p^n}$, and as a result we
obtain Corollary \ref{cor_meta_abel}. In subsection \ref{rel_free}, we
prove Corollary \ref{col_free_nil_p}, and in subsection \ref
{auto_free_nil} we study the structure of the group of tame automorphisms
of the free nilpotent group and calculate explicitly its Kazhdan
constants.

\section{Proof of Theorem \ref{thm_gen_ext}}\label{abel_ext}

% Let $(\pi,V)$ be a unitary representation of $\Gamma$ which don't contains the trivial
%representation. Decompose it into irreducible components, and let
Let $(\rho, \mathcal{H} )$ be a non-trivial unitary irreducible
representation of $\Gamma$ and let $$W_0=\{w \in \mathcal{H} \mid
\rho(a)w=w \quad \forall a \in A\}$$ be the space of $A$-invariant
vectors. Since $A \lhd \Gamma$ , $W_0$ is $\Gamma$ invariant. As $\rho$
is irreducible, we have either $W_0=\mathcal{H}$ or $W_0=\{0\}$.

If $W_0=\mathcal{H}$ then $\rho$ gives rise to a representation
$\bar{\rho}$ of $H=\Gamma/A$. It is clear that $\bar{\rho}$ has no
invariant vectors and thus for every $v \in \mathcal{H}$ with $\|v\|=1$,
we get
$$\frac{1}{|S|+|B|}\sum \limits_{ s \in S}\|\rho(s)v-v \|^2=\frac{1}{|S|+|B|}\sum \limits_{ s \in S}\|\bar{\rho}(s)v-v \|^2 >
 \frac{{\ep_H}|S|}{|S|+|B|}.$$

We assume from now on that $W_0=\{0\}$. The restriction of $\rho$
to $A$ decomposes into irreducible representations:
$$\rho=\oplus_{i=1}^n {m_i \rho_i}$$ where each $\rho_i$ is one
dimensional and the $m_i$ are non negative integers. For each $i
\in \{1,...,n\}$, we denote the isotopic component $m_i \rho_i$ by
$\chi_i$ and we write $X=\{\chi_i\}_{i=1}^n$. Then $X$ is in one
to one correspondence with a subset of the set of irreducible
characters of $A$.  By the irreducibility of $\rho$, $H$ acts
transitively on $X$ (by conjugation). This action gives us a
representation $\Phi$ of $H$ on $\ell _2(X)$ (the vector space of
all complex valued functions on $X$ endowed with the standard
inner product) defined by $(\Phi(h)F)(\chi)=F(\chi^h)$.

Now assume that for some ${\ep}$ which will be determined later, there
exists a unit vector $v \in \mathcal{H}$, s.t
$$\frac{1}{|S|+|B|} \sum
\limits_{g \in \tilde{S} \cup B}\| \rho (g)v-v\|^2 \leq {\ep}.$$
Decompose  $v$ into its isotopic components: $v=\sum \limits_{\chi
\in X} v_{\chi}$. We define $F=F_{v} \in \ell _2 (X)$ by $F(\chi
)=\|v_{\chi }\|$. Note that $F$ is a unit vector.

The following lemma is a general property of Kazhdan groups.

\blem \label{ClInv} Let $\Gamma$ be a group generated by a set $S$
with a Kazhdan constant ${\ep}>0$ with respect to $S$. Let $0 <
\delta <{\ep}$ be given. Let $(\rho,\mathcal{H})$ be a unitary
representation of $\Gamma$ having a unit vector $v \in
\mathcal{H}$, satisfying $\frac{1}{|S|}\sum \limits_{s \in
S}\|\rho(s)v-v\|^2 \leq \delta.$ Then there exists an invariant
vector $v_0$ satisfying $\|v-v_0\|^2<\frac{\delta}{{\ep}}$.
 \elem

\begin{proof}

 Decompose $\rho$ into its trivial and non-trivial components,
 $\rho=\sigma_0+\sigma_1$, and accordingly decompose $v=v_0+v_1$.

Now $\sigma_1$ has no invariant vectors, therefore

\begin{equation}\label{eqt1}
\delta \geq  \frac{1}{|S|}\sum \limits_{s \in S} \|\rho(s)v-v\|^2=
\frac{1}{|S|}\sum \limits_{s \in S}\|\sigma_1(s)v_1-v_1\|^2 >
{\ep}\|v_1\|^2.
\end{equation}

This implies that $\|v_1\|^2<\frac{\delta}{{\ep}}$,
 hence we get $\|v-v_0\|^2<\frac{\delta}{{\ep}}$ as required.

\end{proof}

Here is a technical lemma which we will use latter in the proof.

 \blem
\label{tech} Let $\vec{a}=(a_1,...,a_n)\in \mathbb{R}^n$ be a unit
vector and $\vec{b}=(b_1,...,b_n) \in \mathbb{R}^n$. Assume that
there is some $K\in \mathbb{R}$ such that $|b_i| \leq K/\sqrt{n}$
. Then

\[
\sumlim_{i=1}^n|a_i^2-b_i^2| \leq \sumlim_{i=1}^n(a_i-b_i)^2 +
2K\sqrt{\sumlim_{i=1}^n(a_i-b_i)^2}.
\]

\elem

\begin{proof}

$$
\sumlim_{i=1}^n\left(|a_i^2-b_i^2|-(a_i-b_i)^2\right)=\\
 \sumlim_{\stackrel{i=1}{a_i \geq
b_i}}^n\left((a_i^2-b_i^2)-(a_i-b_i)^2\right)+\sumlim_{\stackrel{i=1}{b_i>
 a_i}}^n\left((b_i^2-a_i^2)-(a_i-b_i)^2\right)=$$
$$\sumlim_{\stackrel{i=1}{a_i \geq
b_i}}^n2b_i(a_i-b_i)+\sumlim_{\stackrel{i=1}{b_i >
a_i}}^n2a_i(b_i-a_i)  \leq \frac{2K}{\sqrt{n}}\sumlim_{\stackrel{i=1}{a_i
\geq b_i}}^n(a_i-b_i)+\frac{2K}{\sqrt{n}}
\sumlim_{\stackrel{i=1}{b_i >
a_i}}^n(b_i-a_i)$$
$$=\frac{2K}{\sqrt{n}}\sumlim_{i=1}^n|a_i-b_i|
 \leq 2K\sqrt{\sumlim_{i=1}^n(a_i-b_i)^2}.
$$

The last inequality is by the Cauchy Schwarz inequality.
\end{proof}

 \blem

$\frac{1}{|S|}\sum \limits_{s \in S}\|\Phi(s)F-F\|^2 \leq
\frac{{\ep}(|S|+|B|)}{|S|}.$ \elem
\begin{proof}
From the definition of $\Phi$, we have:

$$\frac{1}{|S|+|B|}\sum \limits_{s\in S}\|\Phi(s)F-F \|^2= \frac{1}{|S|+|B|}\sum \limits_{s\in S} {\sum_{\chi  \in X} |F(\chi^s)-F(\chi)|^2}
=$$
$$\frac{1}{|S|+|B|}\sum \limits_{s\in S} {\sum_{\chi \in X} |
(\|v_{\chi^s} \| - \|v_\chi \| )|^2} =\frac{1}{|S|+|B|}\sum
\limits_{s\in S} {\sum_{\chi  \in X} | (\|v_{\chi^s} \| - \|
\rho(s)v_\chi \|) |^2}$$ $$ \leq \frac{1}{|S|+|B|}\sum \limits_{s\in S}
{\sum_{\chi  \in X} \| \rho(s)v_\chi-v_{\chi^s}\|^2} \leq {\ep} .$$

\end{proof}

Denote by $c=\frac{|B|}{|S|}$, now by Lemma \ref{ClInv} we obtain that
there exists an invariant (i.e. constant) function $G$ such that $\| F-G
\|^2 < \frac{(1+c)\cdot{\ep}}{{\ep}_H}$. Thus
$$\|G\|=K<\|F\|+\sqrt{\frac{(1+c)\cdot{\ep}}{{\ep}_H}}=1+\sqrt{\frac{(1+c)\cdot{\ep}}{{\ep}_H}}$$
and so $G(\chi )= K/\sqrt{|X|}$. We choose $u \in \mathcal{H}$
such that for each $\chi \in X$, $$G(\chi)=\|u_{\chi}\|.$$

Since $W_0=\{ 0\}$, the following holds:
$$ {\ep}_{A} < \frac{1}{|B||H|} \sum \limits_{b \in B} \sum \limits_{h \in H} \| \rho(b^h)v-v \|^2 = \frac{1}{|B||H|} \sum \limits_{b \in B} \sum \limits_{h \in H} \sum_{\chi \in X} \| \chi(b^h)v_\chi-v_\chi\|^2= $$

$$ \frac{1}{|B||H|} \sum \limits_{b \in B} \sum \limits_{h \in H} \sum_{\chi \in X} |\chi^h(b)-1|^2\|v_\chi\|^2 = \frac{1}{|B||H|} \sum \limits_{b \in B} \sum \limits_{h \in H} \sum_{\chi \in X} |\chi(b)-1|^2 F(\chi^{h^{-1}}) ^2.$$

Now,

$$ |\frac{1}{|B||H|} \sum \limits_{b \in B} \sum \limits_{h \in H} \sum_{\chi \in X} | \chi(b)-1|^2 F(\chi^{h^{-1}})^2 - \frac{1}{|B||H|} \sum \limits_{b \in B} \sum \limits_{h \in H}\sum_{\chi \in X} | \chi(b)-1|^2 G(\chi )^2 | \leq $$

$$ \frac{1}{|B||H|} \sum \limits_{b \in B} \sum \limits_{h \in H} \sum_{\chi \in X} | \chi(b)-1|^2 | F(\chi^{h^{-1}})^2 - G(\chi )^2 | \leq$$

$$4\cdot \frac{1}{|B||H|} \sum \limits_{b \in B} \sum \limits_{h \in H} \sum_{\chi \in X} | F(\chi^{h^{-1}})^2-G(\chi )^2|=
4\cdot \frac{1}{|B||H|} \sum \limits_{b \in B} \sum \limits_{h \in H} \sum_{\chi \in X} | F(\chi^{h})^2-G(\chi )^2|.$$

Now, substituting $a=F$ and $b=G$ in Lemma \ref{tech}, we get:
$$ 4\cdot \sum_{\chi \in X} | F(\chi^h)^2-G(\chi )^2|
\leq 4 \cdot \left(\sumlim_{\chi \in X
}(F(\chi^{h})-G(\chi ))^2+2K\sqrt{\sumlim_{\chi
\in X }(F(\chi^h)-G(\chi ))^2}\right)=$$

$$4 \cdot \left(\sumlim_{\chi \in X }(F(\chi )-G(\chi ))^2+2K\sqrt{\sumlim_{\chi \in X }(F(\chi )-G(\chi
))^2}\right)\leq 4 \cdot \left(\|F-G\|^2+2K\sqrt{\|F-G\|^2}\right) $$

$$<4 \cdot \left(\frac{(1+c)\cdot{\ep}}{{\ep}_{H}} +2
(1+\sqrt{\frac{(1+c)\cdot{\ep}}{{\ep}_{H}}})\sqrt{\frac{(1+c)\cdot{\ep}}{{\ep}_{H}}}\right)=
12\cdot \frac{(1+c)\cdot{\ep}}{{\ep}_{H}}+8\cdot\sqrt{\frac{(1+c)\cdot{\ep}}{{\ep}_{H}}}$$

Now

$$ {\ep}_A < \frac{1}{|B||H|} \sum \limits_{b \in B} \sum \limits_{h \in H} \|\rho(b^h)v-v\|^2 = \frac{1}{|B||H|} \sum \limits_{b \in B} \sum \limits_{h \in H}\sum_{\chi \in X} |\chi(b)-1|^2\|v_{\chi^{h^{-1}}}\|^2 \leq $$
$$ 12\cdot \frac{(1+c)\cdot{\ep}}{{\ep}_H}+8\sqrt{\cdot\frac{(1+c)^2\cdot{\ep}}{{\ep}_H}} +
\frac{1}{|B||H|} \sum \limits_{b \in B} \sum \limits_{h \in H}  \sum_{\chi \in X} |\chi(b)-1|^2G(\chi )^2 \leq $$

$$12\cdot
\frac{(1+c)\cdot{\ep}}{{\ep}_H}+8\sqrt{\cdot\frac{(1+c)\cdot{\ep}}{{\ep}_H}}+12\cdot
\frac{(1+c)\cdot{\ep}}{{\ep}_H}+8\sqrt{\cdot\frac{(1+c)\cdot{\ep}}{{\ep}_H}}
+\frac{1}{|B||H|} \sum \limits_{b \in B} \sum \limits_{h \in H} \sum_{\chi \in X} |\chi(b)-1|^2\|v_{\chi}\|^2 =$$

$$ 24\cdot
\frac{(1+c)\cdot{\ep}}{{\ep}_H}+16\sqrt{\cdot\frac{(1+c)\cdot{\ep}}{{\ep}_H}}
+\frac{1}{|B||H|} \sum \limits_{b \in B} \sum \limits_{h \in H}
\sum_{\chi \in X} |\chi(b)-1|^2\|v_{\chi}\|^2 =$$
$$  24\cdot
\frac{(1+c)\cdot{\ep}}{{\ep}_H}+16\cdot\frac{(1+c)\cdot{\ep}}{{\ep}_H}
+\frac{1}{|B|} \sum \limits_{b \in B} \sum_{\chi \in X}
|\chi(b)-1|^2\|v_{\chi}\|^2 < $$

 $$24\cdot
\frac{(1+c)\cdot{\ep}}{{\ep}_H}+16\sqrt{\cdot\frac{(1+c)\cdot{\ep}}{{\ep}_H}}+
(1+\frac{1}{c}){\ep}.$$

If we choose $${\ep} \leq \frac{{\ep}_H
{\ep}_A}{512\cdot(1+c+\frac{1}{c}) } ,$$ we get a contradiction (the last
inequality follows from the fact that ${\ep},{\ep}_H \leq 1$) and this
finishes the proof of Theorem \ref{thm_gen_ext}.

 \qed

\section{Applications}

\subsection{Solvable Groups} \label{solv_grp}

Lubotzky and Weiss (\cite{LW}, Corollary 3.3) proved that Kazhdan
constants of an infinite family of solvable groups with bounded derived
length, with respect to any bounded generating set are going to zero. The general
way to give a positive Kazhdan constant for an infinite family of solvable groups is by using  Alon-Roichman result which require a generating set of logarithmic size of the group.

The first improvement for some solvable groups is due to Meshulam and
Wigderson \cite{MW} (based on \cite{RVW,ALW}). They constructs a family of
solvable groups of the form $G_{i+1}=G_i \ltimes \mathbb{F}_{p_i}[G_i]$ where
$G_i$ is a finite group, $p_i$ is a prime and $G_{i+1}$ is the natural
semi-direct product.

\bthm  (\cite{MW}, Theorem 1.7) There exists a group $G_1$ and a sequence
of primes $p_i$, such that for any $0 <\ep<1$ we have $$g_{\ep}(G_n)\leq
\log^{(n-\log^*n)}|G_n|,$$ where $\log^k$ denotes the $k$ times iterated
logarithm function. \ethm

We would like to mention a related work by Rozenman, Shalev and Wigderson \cite{RSW} (which is also based on \cite{RVW,ALW}).
In this work they gave an iterative construction of infinite family of finite (non-solvable) groups $\{G_n\}_{n=1}^\infty$ and they showed that 
for any $n\in \mathbb{N}$, $g_\ep(G_n) = O(1)$ for some $\ep>0$.

 Influenced by Meshulam and Wigderson methods, we will
give in this section a construction  of a family of meta-abelian groups
$\Gamma_p$ (split and non-split), with a substantial improvement
for $g_\ep(\Gamma_p)$ (compared to Alon-Roichman).

\subsection{Metabelian Groups}

\subsubsection{Group Cohomology and Non-Splitting Extension}
\begin{defn} If $\pi:\Gamma \rightarrow H$ is surjective homomorphism of groups, then a
\textit{lifting} of $x\in H$ is an element $l(x)\in \Gamma$ with
$\pi(l(x))=x$.
\end{defn}

Let $A$ be a group, $H$ a group with a homomorphism
$$\theta:H\rightarrow Aut(A).$$
Let $H^2(H,A,\theta)$ be the second
cohomology group. As is well known (see \cite{Rot}, ch. 7),
$H^2(H,A,\theta)$ parameterizes groups extension of $A$ by $H$ with respect to
$\theta$.

Let $p$ be a prime, $1 \leq n \in \mathbb{N}$ and let
$C_{p^n}$ be the cyclic group of order $p^n$.
In the next lemma we will prove the existence of non-splitting extension of a $p$-group by some cyclic groups.
Next, by using this lemma we will construct a families of metabelian groups with non-trivial second cohomology.

\blem Let $p$ be a prime, $2\leq n\in \mathbb{N}$ and $F$ a non-trivial
finite $p$-group. Assume that the cyclic group $C_{p^n}$ has non-trivial action $\theta$ on $F$ and
for some $k<n$, $p^k \cdot C_{p^n}$ acts trivially on $F$, and $F$ has
exponent dividing $p^{n-k}$ (i.e. $x^{p^{n-k}}=1$ for all $x \in F$). Then $$H^2(C_{p^n},F,\theta)\neq 0.$$ \elem

\begin{proof} Let ${\Gamma}= C_{p^n}\ltimes F$ which is a finite $p$-group, $F$ is a non-trivial normal subgroup
so has non-trivial intersection with the center. So there exists a
subgroup $Z$ of $F$, cyclic of order $p$, which is central in $C_{p^n}\ltimes F$. Now
the cyclic group $C_{p^{n+1}}$ has a unique subgroup of order p which we denote by $Z'$.
Let $\tilde{\Gamma}=C_{p^{n+1}}\ltimes F$. Consider a
"diagonal subgroup"  $K$ in $Z\times Z'$. Then it is easy to check that $F$ is a normal subgroup of
$\tilde{\Gamma}/K$. So $\tilde{\Gamma}$ is an extension of $F$ by $C_{p^n}$,
We claim that the generator of $C_{p^n}$ cannot be lifted to an element of
order $p^n$. Indeed, if $(c,f)$ is a lifting, then taking the power $p^k$ we
get an element $(c^{p^k},f')$ and $c^{p^k}$ centralizes $f'$, so taking the
power $p^{n-k}$ of the latter (i.e. the power $p^n$ of the initial element)
and using the exponent, we obtain $(c^{p^n},1)$; this element of $C_{p^{n+1}}\ltimes
F$ remains non-trivial after taking the quotient by $K$.

\end{proof}

 Let $\mathbb{F}_p$ be the field of order $p$ and
$\mathbb{F}_p[C_p]$ the group algebra. Every element $f \in
\mathbb{F}_p[C_p]$ can be represented as $f=\sum \limits_{x \in
C_p}f(x)x$.

Let $\theta:C_{p^n}\rightarrow Aut(\mathbb{F}_p[C_p])$ be the left translation
(mod $p$) homomorphism.
As a consequence of last lemma (take $k=1$) we get the following result:

 \blem \label{non_splitting} For every $n \geq 2$
$$H^2(C_{p^n},\mathbb{F}_p[C_p],\theta) \neq 0,$$
hence, there is a non-split extension of $A=\mathbb{F}_p[C_p]$ by $H=C_{p^n}$.
\elem

\subsubsection{Proof of Corollary \ref{cor_meta_abel}:}

We begin by proving that for any $0<\ep<1$, there exist a subset $B
\subseteq \mathbb{F}_p[C_p]$ such that the average Kazhdan constant of $\mathbb{F}_p[C_p]$ with respect to
$\tilde{B}=B^{C_p}$ is greater than $\ep$.

  Given $f,g \in \mathbb{F}_p[C_p]$, define the product
$f\cdot g = \sum \limits_{x\in C_p}f(x)g(x) \in \mathbb{F}_p$.
Let $e_p(\alpha)=\exp(\frac{2\pi \alpha i}{p})$. A multiset $ \tilde{B}
\subset \mathbb{F}_p[C_p]$ is called $\delta$-balanced if for all $0 \neq f \in
\mathbb{F}_p[C_p]$
$$|\sum \limits_{h \in \tilde{B}}e_p(f \cdot h)| \leq (1-\delta)|\tilde{B}|.$$

If $\tilde{B}$ is $\delta$-balanced then
$$1-(\frac{1}{|\tilde{B}\cup\tilde{B}^{-1}|}|\sum \limits_{h \in \tilde{B}\cup\tilde{B}^{-1} } e_p(f\cdot h)|)=
1-\frac{1}{2|\tilde{B}|}|\sum \limits_{h \in \tilde{B} } e_p(f\cdot
h)+\sum \limits_{h \in \tilde{B}^{-1} } e_p(f\cdot h)| \geq$$

$$1-\frac{1}{2|\tilde{B}|}(|\sum \limits_{h \in \tilde{B} } e_p(f\cdot
h)|+|\sum \limits_{h \in \tilde{B}^{-1} } e_p(f\cdot h)|) = 1-\frac{1}{|\tilde{B}|}|\sum \limits_{h \in \tilde{B} } e_p(f\cdot
h)| \geq \delta. $$

Since $\mathbb{F}_p[C_p]$ is an abelian group, it is easy to check that the
spectral gap of $\mathbb{F}_p[C_p]$ with respect to $\tilde{B}\cup\tilde{B}^{-1}$ is:
$$\min \limits_{\chi \neq 1}|\frac{1}{|\tilde{B}\cup\tilde{B}^{-1}|}\sum \limits_{h \in \tilde{B}\cup\tilde{B}^{-1} }\chi(h)-1 |.$$

So, from the discussion in subsection \ref{spectral_gap}, we obtain that
if $\tilde{B}$ is $\delta$-balanced, then the average Kazhdan constant of
$A$ with respect to $\tilde{B}\cup\tilde{B}^{-1}$ is at least $2\delta$.

 For $f\in \mathbb{F}_p[C_p]$, $s \in \mathbb{N}$ and $\delta >0$ let
 $$B_\delta(f)=\left\{(h_1,...,h_s)\in
\mathbb{F}_p[C_p]^s:\left|\frac{1}{sp}\sum_{i=1}^s\sum_{\sigma \in C_p}e_p(\sigma
h_i\cdot f)\right| > 1- \delta\right\}.$$

\noindent For $f\in \mathbb{F}_p[C_p]$, define a linear map $T_f:\mathbb{F}_p[C_p]\rightarrow
\mathbb{F}_p[C_p]$ by $T_f(h)=hf$.
Denote by $C_pf$ the orbit of $f$ under the action of $C_p$. It is clear
that $$\dim (Span \{C_pf\})=\textnormal{rank}T_f.$$
 Let $$V_r(\mathbb{F}_p)=\{f\in \mathbb{F}_p[C_p]: \textnormal{rank}T_f = r\}.$$

Let $C_{p^n}$  be the cyclic group of order $p^n$ for some $n \in \mathbb{N}$. We defined above the following representation of $C_{p^n}$: $$\theta:C_{p^n} \rightarrow Aut(\mathbb{F}_p[C_p]).$$
Let $a$ be a generator of $C_{p^n}$, set $\alpha=\theta(a)$, then $\alpha$ satisfy $\alpha^{p^n}=1$.
 Thus $\alpha$
is a root of the polynomial $x^{p^n}-1=(x-1)^{p^n}$. So the minimum
polynomial of $\alpha$ is $(x-1)^k$ for some $ k\leq p^n$. Hence, $1$ is
the only eigenvalue of $\alpha$, so the trivial representation (i.e. this is the subspace of constant vectors)
is the only irreducible representation of
$C_{p^n}$.

Now it is immediate to check that $\mathbb{F}_p[C_p]$ is an indecomposable
representation of dimension $p$. Furthermore, consider the Jordan
canonical form of $\alpha$, we see that a subspace of $\mathbb{F}_p[C_p]$ is
indecomposable if and only if it correspond to a Jordan block which may
be of size $1\leq s \leq p$. Therefore, if one represent $\alpha$ in it
Jordan canonical from, then we get that for every $1\leq k\leq p$ the
number of elements in $\mathbb{F}_p[C_p]$ of rank $k$, is $(p-1)p^{k-1}$.

Meshulam and Wigderson (see Proposition 3.2 in \cite{MW}) proved the
following:
 \bpr \label{good_orbit} If $\textnormal{rank}T_f=r$ then
$$Prob(B_\delta(f)) \leq 8\exp(\frac{-(1-2\delta)^2rs}{4}).$$\epr

The following proposition is a variant of Theorem $1.2$ in \cite{MW}.
 \bpr \label{pr_orbits}For any $0<\delta <\frac{1}{2}$ there exist
 $s=O\left(\frac{1}{(1-2\delta)^2}\cdot \ln p\right)$ elements $h_1,...,h_s\in
 \mathbb{F}_p[C_p]$ such that the multiset
 $\tilde{B}=\bigcup_{i=1}^sC_{p}h_i\subset \mathbb{F}_p[C_p]$ is
 $\delta$-balanced and the average Kazhdan constant of $\mathbb{F}_p[C_p]$ with respect to $\tilde{B}\cup\tilde{B}^{-1}$ is $2\delta$.
 \epr

\begin{proof}

$$Pr\left(\bigcup \limits_{0 \neq f \in \mathbb{F}_p[C_p]} B_\delta(f)\right) \leq \sum \limits_{0 \neq f \in \mathbb{F}_p[C_p]}
B_\delta(f)\leq$$
$$8 \sum \limits_{r\geq 1}
|V_r(\mathbb{F}_p)|\exp\left(\frac{-(1-2\delta)^2rs}{4}\right) \leq$$

$$8 \sum \limits_{r=1}^p p^r\exp\left(\frac{-(1-2\delta)^2rs}{4}\right)=
 8\sum \limits_{r=1}^p(p^4\exp(-(1-2\delta)^2 s)^\frac{r}{4}=$$

 $$8 [p^4\exp(-(1-2\delta)^2
s)]^{\frac{1}{4}}\frac{1-[p^4\exp(-(1-2\delta)^2
s]^{\frac{p}{4}}}{1-[p^4\exp(-(1-2\delta)^2 s]^{\frac{1}{4}}} <
1$$

For any $c> 1$ and for $p$ large enough, if one choose
$s=\frac{4}{(1-2\delta)^2}(c\ln p)$, then the probability that
$h_1,...,h_s$ is not $\delta$-balanced is strictly less than $1$ and the
result follow.
\end{proof}

The last proposition still holds if we replace the group $C_p$ by
$C_{p^n}$ with the action by left translation mod $p$ on
 $\mathbb{F}_p[C_p]$. This follow easily from the fact that for every $f\in \mathbb{F}_p[C_p]$,
every element in the multiset $C_pf$ appears $p^{n-1}$ times in the
multiset $C_{p^n}f$.

Now we are ready to complete the proof of Corollary \ref{cor_meta_abel}.
By \cite{AR} for any $0<{\ep}_1<1$, there is a generating set $S$ for
$C_{p^n}$ of size $O_{{\ep}_1}(\log|C_{p^n}|)$ such that the Kazhdan
constant of $C_{p^n}$ with respect to $S$ is bigger than ${\ep}_1$. Let $\Gamma$
be any extension of $\mathbb{F}_p[C_p]$ by $C_{p^n}$ with respect to to the action of left
translation. Let $\tilde{S}$ be any lifting of $S$ to $\Gamma$, from
Proposition \ref{pr_orbits} (and the comment afterward), we get that for
any $0< {\ep}_2 < 1$ there exist a set $B \subset \mathbb{F}_p[C_p]$ such that the
Kazhdan constant of $\mathbb{F}_p[C_p]$ with respect to $\tilde{B}=B^{C_{p^n}}$ is greater
then ${\ep}_2$. Now from Theorem \ref {thm_gen_ext} the proof is
complete.
 \qed

\subsection{$p$-Groups} \label{rel_free}

Fix a prime $p$, for any group $\Gamma$, define $\Gamma^p$ to be the
subgroup generated by $\{g^p:g\in \Gamma\}$. The lower $p$-series (also
called the lower central $p$-series or the lower exponent-$p$ series) of
$\Gamma$ is the descending series
$$\Gamma=\varphi_1(\Gamma) \geq \varphi_2(\Gamma)\geq \cdot\cdot\cdot
\geq \varphi_n(\Gamma)\geq \cdot\cdot\cdot ,$$ defined by
$$\varphi_{i+1}(\Gamma)=\varphi_i(\Gamma)^p[\varphi_i(\Gamma),\Gamma]$$
for $i\geq 1$. The group $\Gamma$ is said to have class $n$, if
$\varphi_n(\Gamma)$ is the last non-identity element of the lower
$p$-series.

This subsection uses results and properties of the lower $p$-series.
Proofs of the results can be found in Huppert and Blackburn \cite{HB}
(ch. VIII) and in Blackburn, et el. \cite{BNV} (ch. 4). On can
check the following:
 \bpr \label{p_series_center}  For all positive integer
 $i$,
$\varphi_{i+1}(\Gamma)$ is the smallest normal subgroup of
     $\Gamma$ lying in $\varphi_i(\Gamma)$, such that
     $\varphi_i(\Gamma)/\varphi_{i+1}(\Gamma)$ is an elementary
     abelian $p$ group and is central in
     $\Gamma/\varphi_{i+1}(\Gamma)$.

 \epr

Let $F_k$ be the free group on $k$ generators. Let
$\Gamma_{k,n}=F_k/\varphi_{n+1}(F_k)$. This group is a $p$-group,
called the relatively free group on $k$ generators in the variety
of $p$-groups with a lower-$p$ series of class $n$ (for an
introduction to variety of a groups see ch. 1 in \cite{Ne}). The
dimensions of $\varphi_i(F_k)/\varphi_{i+1}(F_k)$ can be
calculated explicitly (see \cite{HB}), in particular for the case
of $\Gamma_{k,2}=F_k/\varphi_3(F_k)$, we have the following lemma:

\blem (\cite{BNV}, Lemma 4.2) The Frattini subgroup $\Phi(\Gamma_{k,2})$
of $\Gamma_{k,2}$ is of order $p^{\frac{1}{2}k(k+1)}$ and index $p^k$.
\elem

As a consequence of the last lemma, one can show that the
order of the commutator subgroup $\Gamma_{k,2}'$ is
$p^{\frac{1}{2}k(k-1)}$.

\subsubsection {Proof of Corollary \ref{col_free_nil_p}:}
 Let $\Gamma_{k,2}$ be the
relatively free group in the variety of $p$-groups with a lower-$p$
series of class. The Frattini subgroup of $\Gamma_{k,2}$ is
$$\Phi(\Gamma_{k,2})={\Gamma_{k,2}}^p[\Gamma_{k,2},\Gamma_{k,2}]=\varphi_2(\Gamma_{k,2}).$$
Hence if $\tilde{S}\subset \Gamma_{k,2}$ is such that its image in
$\Gamma_{k,2}/\varphi_2(\Gamma_{k,2})$ generates
$\Gamma_{k,2}/\varphi_2(\Gamma_{k,2})$, then $\tilde{S}$ generates
$\Gamma_{k,2}$.

From Proposition \ref{p_series_center}, $\varphi_2(\Gamma_{k,2})$ is
central in $\Gamma_{k,2}$, so $\Gamma_{k,2}$ is a central extension of
$\Gamma_{k,2}'=[\Gamma_{k,2},\Gamma_{k,2}]$ by
$\Gamma_{k,2}/\Gamma_{k,2}'$.

By \cite{AR} for any $0<\ep <1$ there is a generating set $S$ for
$\Gamma_{k,2}/\Gamma_{k,2}'$ of size
$O_{{\ep}}(\log|\Gamma_{k,2}/\Gamma_{k,2}'|)=O_{{\ep}}(\log|p^{2k}|)$
such that the average Kazhdan constant of $\Gamma_{k,2}/\Gamma_{k,2}'$ with respect to $S$
is bigger than ${\ep}$. Let $\tilde{S}$ be any lifting of $S$ to
$\Gamma_{k,2}$. We claim that the average Kazhdan constant of
$\Gamma_{k,2}$ with respect to $\tilde{S}$ is bigger than ${{\ep}} /4$. Indeed,
one dimensional representations are trivial on $\Gamma_{k,2}'$ (since it
is the commutator subgroup) and hence factor through
$\Gamma_{k,2}/\Gamma_{k,2}'$, and for higher dimensional representations
we use theorem \ref{thm_cen_ext}. We have shown that $\Gamma_{k,2}$ has a
generating set of size $$O_{{\ep}}(\log
|\Gamma_{k,2}/\Gamma_{k,2}'|)=O_{{\ep}}(\sqrt{\log |\Gamma_{k,2}|}).$$
One can repeat the above methods and improve $g_\ep(\Gamma_{k,n})$  for
the relatively free group on $k$ generators in the variety of $p$-groups
with a lower-$p$ series of class $n>2$.

\qed

\subsection{Kazhdan Constants for the Tame Automorphism Groups of Free Nilpotent Groups} \label{auto_free_nil}

Let $F_k$ be the free group on $k$ generators $x_1,...,x_k$. For any $1
\leq i\neq j \leq k$, let $R_{i,j}^{\pm},L_{i,j}^{\pm}$ be the following
automorphisms of $F_k$:
$$R_{i,j}^{\pm}(x_i)=x_ix_j^{\pm1} \mbox{ and } R_{i,j}^{\pm}(x_l)=x_l \mbox{ if }
l\neq i.$$
$$L_{i,j}^{\pm}(x_i)=x_j^{\pm1}x_i \mbox{ and } L_{i,j}^{\pm}(x_l)=x_l \mbox{ if }
l\neq i.$$
 These automorphisms are called Nielsen moves. Let $S$ be
the set of Nielsen moves of the free group on $k$ generators, it is easy
to check that $|S|=4k(k-1)$.

\begin{defn} Let $F_k$ be the free group on $k$ generators and
$W=\gamma_{c+1}(F_k)$ be the $(c+1)$-th term of the lower central series
of $F_k$. Then $F_k/W$ is called the free nilpotent group of class $c$
and will be denoted by $F_k(c)$. \end{defn}

Let $U$ be a characteristic subgroup of $F_k$, then it is clear
that every automorphism of $F_k$ induces
 an automorphism on $F_k/U$. This fact bring us to the following definition.

\begin{defn} An automorphism of the free nilpotent group $F_k(c)$ is said to be
\emph{tame} if it is in the image of the natural projection map
$Aut(F_k)\rightarrow Aut(F_k(c))$. The subgroup of all the tame
automorphisms is denoted by $A_k(c)$. \end{defn}

The group $A_k(c)$ has the following short exact sequence:
$$1 \rightarrow K_1 \rightarrow A_k(c)\rightarrow SL_k(\mathbb{Z}) \rightarrow 1$$
where $K_1$ is a nilpotent group of class $c-1$ (see \cite{And} for
details).

 Lubotzky and Pak showed in \cite{LP} that for $k \geq
3$, the group $A_k(c)$ is a lattice in some Lie group which has
Kazhdan property (T), hence  it has property (T). Their approach
does not give an estimate for the Kazhdan constant. Mallahi
\cite{MAK} following Burger, Shalom and Kassabov (see
\cite{Bur,Sh1,Kas1}), used the spectral measure on the dual group
$\hat{\mathbb{Z}}^m$ correspond to $SL_k(\mathbb{Z})$ and prove
that the group $A_k(c)$ has relative property (T) with respect to
the center of $K_1$. Before we state Mallahi results in a more
general setting (which follows easily from Mallahi's proof), we
give some notation. Let $k,c \in \mathbb{N}$, we define
$$\delta(k,c)=\frac{2\sqrt{c}}{\sqrt{(84\sqrt{k}+1920)(4k)^c}}.$$

 \bthm (Theorem 2.5.6 in \cite {MAK}) \label{thm_mak_2} Let $k\geq 3$
be an integer, $A_k(c)$ be the group of tame automorphisms of the free
nilpotent group of class $c$ and $S$ the set of Nielsen transformations.
Let $A$ be the center of $A_k(c)$ and let $Z$ be the center of $K_1/A$.
Let $(\rho,\mathcal{H})$ be a unitary representation of $A_k(c)/A$
 with $(S,\delta(k,c))$-invariant vector.
 %where
% $$\delta(k,c)=\frac{2\sqrt{c}}{\sqrt{(84\sqrt{k}+1920)(4k)^c}}.$$
Then $\mathcal{H}$ contains a non-zero $Z$-invariant vector. \ethm

From the last theorem, Mallahi deduce the following:
 \bthm (Theorem 2.5.7 in \cite {MAK}) \label{thm_mak_1}
 Let $k\geq 3$ be an integer, $A_k(c)$ be the
 group of tame automorphisms of the free nilpotent group of class $c$ and $S$ the set of Nielsen
transformations. If $c < 2k+1$, then the Kazhdan constant
$\epsilon(k,c)$ of $A_k(c)$ with respect to $S$ is greater than or
equal to
$$ \frac{\sqrt{c}}{\sqrt{(84\sqrt{k}+1920)(4k)^c}}.$$ \ethm

 The limitation for $c<2k+1$ is based on the following fact due to
Formanek (Theorem 6 in \cite{For}):

\bthm  \label{thm_For} Let $F_k(C)$ be the free nilpotent group of rank
$k$ and class $c$ where $k,c \geq 2$. Then the automorphism group has a
non-trivial center iff $c=2kl+1$ for $l\geq 1$. \ethm

The goal of this section is to extend Mallahi result and calculate
Kazhdan constant of $A_k(c)$ for any $k \geq 3$ and $c\geq 1$. Using the
same notations as above we prove:

\bthm  Let $k\geq 3$ be an integer, $A_k(c)$ be the
 group of the tame automorphisms of the free nilpotent group of class $c$ and $S$ the set of Nielsen
transformations, then the Kazhdan constant $\ep(k,c)$  of $A_k(c)$ with respect to  $S$ satisfies:

$$\epsilon(k,c) \geq \left\{ \begin{array}{ll}

        \frac{\sqrt{c}}{\sqrt{(84\sqrt{k}+1920)4^{c+3}k^{c+4}}} &
        \mbox{if $c = 2kl+1 $ for some $l\in\mathbb{N}$};\\

        \frac{\sqrt{c-1}}{\sqrt{(84\sqrt{k}+1920)4^{c+3}k^{c+3}}} &
        \mbox{if $c = 2kl+2 $ for some $l\in\mathbb{N}$};\\

        \frac{\sqrt{c-2}}{\sqrt{(84\sqrt{k}+1920)4^{c+3}k^{c+2}}} &
        \mbox{if $c = 2kl+3 $  for some $l\in\mathbb{N}$};\\

        \frac{\sqrt{c-3}}{\sqrt{(84\sqrt{k}+1920)4^{c+3}k^{c+1}}} &
        \mbox{if $c = 2kl+4 $  for some $l\in\mathbb{N}$} \\
        \frac{\sqrt{c-4}}{\sqrt{(84\sqrt{k}+1920)4^{c+3}k^c}} &
        \mbox{otherwise}.\end{array} \right.$$ \ethm

\begin{proof}

For the case $c  \leq 2k$, by Theorem \ref{thm_mak_1} we get that the
Kazhdan constant of $A_k(c)$ with respect to $S$ is
$$\epsilon(k,c) \geq \frac{\sqrt{c}}{\sqrt{(84\sqrt{k}+1920)(4k)^c}}.$$

For the case $c =2k+1$, denote by $A$ the center of $A_k(c)$, by Theorem
\ref{thm_For}, $A$ is not trivial. Hence we have the following exact
sequence:
$$ 1 \rightarrow A \rightarrow A_k(c)\rightarrow A_k(c)/A \rightarrow 1,$$
and respectively
$$ 1 \rightarrow K_1/A \rightarrow A_k(c)/A \rightarrow SL_k(\mathbb{Z}) \rightarrow 1.$$

\blem \label{A_k_Kaz} Let $c=2k+1$ then the Kazhdan constant of
$A_k(c)/A$ is greater than or equal to $\frac{\delta(k,c)}{2}$.
%$$\frac{\sqrt{c}}{\sqrt{(84\sqrt{k}+1920)(4k)^c}}.$$
\elem
\begin{proof}

Let $(\rho,,\mathcal{H})$ be a unitary representation of $A_k(c)/A$ and
assume that there exist a unit vector $v \in \mathcal{H}$ which is
$(\frac{\delta(k,c)}{2},S)$-invariant. Let $Z$ be the center of $K_1/A$,
put
$$\mathcal{H}_0=\{ v \in \mathcal{H}: \rho(z)v=v, \forall z \in Z \}$$

\noindent and $\mathcal{H}_1=\mathcal{H}_0^\perp$ (the orthogonal complement).
Write $v=v_0+v_1$ where $v_i\in\mathcal{H}_i$, by Theorem \ref{thm_mak_2}
we get that $\mathcal{H}_0\neq {0}$.

The subspace $\mathcal{H}_1$ does not have a $Z$-invariant vector,
and therefore it does not have a $(\delta(k,c),S)$-invariant
vector. This implies that $\delta(k,c)\|v_1\| < \max \limits_{s \in
S} \|\rho(s)v_1-v_1\| \leq \frac{\delta(k,c)}{2}$, and we deduce
that $\|v_1\|  < \frac{1}{2}$.

So $\|v_0\|\geq \frac{1}{2}$ and $\mathcal{H}_0$ is $A_k(c)/A$-invariant
because $Z$ is a normal subgroup. Furthermore, $Z=\gamma_{c}(K_1/A)$,
therefore it rise a unitary representation $\overline{\rho}$ of
$$(A_k(c)/A)/Z=A_k(c-1)$$ on $\mathcal{H}_0$.

Let $\tilde{S}$ be the projection of the generators $S$ to
$A_k(c-1)$, then it is clear that $\tilde{S}$ is the set of
Nielsen transformations for $A_k(c-1)$. Now, for any $s \in S$ we
have

$$\|\rho(s)v-v\|^2=\|\overline{\rho}(s)v_0-v_0\|^2 + \|\rho(s)v_1-v_1\|^2 \leq
\frac{\delta(k,c)^2}{4} $$ and therefore
$$\|\overline{\rho}(s)v_0-v_0\|^2 \leq
\frac{\delta(k,c)^2}{4} \leq \delta(k,c)^2\|v_o\|^2 \leq \frac{\delta(k,c-1)^2}{4}\|v_0\|^2.$$

Consequently, the vector $v_0\in \mathcal{H}_0$ is
$(\frac{\delta(k,c-1)}{2},\tilde{S})$-invariant . By Theorem
\ref{thm_mak_1} we obtain that $\mathcal{H}_0$ contains a non-zero
invariant vector of $A_k(c-1)$ and this vector is also an
invariant vector of $A_k(c)/A$.
\end{proof}

The group $A_k(c)$ does not have a non-trivial representation of
dimension one, because $A_k(c)=[A_k(c),A_k(c)]$ (see \cite{Ge})
and $|S|=4k(k-1)$, therefore from Theorem \ref{thm_cen_ext} we
obtain that the Kazhdan constant $\ep(k,c)$ of $A_k(c)$ for
$c=2k+1$ with respect to $S$ is greater than or equal to
$$\frac{\delta(k,c)}{2|S|}\geq \frac{\sqrt{c}}{\sqrt{(84\sqrt{k}+1920)4^{c+3}k^{c+4}}}.$$

For $c >2k+1$, we continue by induction on $c$ and repeating the same arguments as in Lemma \ref{A_k_Kaz}.

\begin{comment}Notice that we should take special care for the cases of $2kl+1 \leq c
\leq 2kl+5$ for any $l\in \mathbb{N}$. For example for $c=2kl+2$, we take $\ep(k,2kl+2)=\frac{\sqrt{c-1}}{\sqrt{(84\sqrt{k}+1920)4^{c+3}k^{c+3}}}$.
So if we have a unitary representation of $A_k(c)$ with a unit vector which is $(S,\ep(k,2kl+2))$-invariant.
Then this representation factor throw a central subgroup $Z$, and we get that there exists a unit vector in a representation of $A_k(2kl+1)$ which is $(S,\ep(k,2kl+1))$-invariant($A_k(2kl+2)/Z=A_k(2kl+1))$).\end{comment}
\end{proof}

\section {Acknowledgments.} This paper is part of the author's
Ph.D. studies under the guidance of Prof. Alex Lubotzky whom I
would like to thank for introducing me the problem and for
fruitful conversations. I am grateful to Nir Avni for his
assistance in proving Theorem \ref{thm_gen_ext}. I would also like to thank Eli Bagno, Martin Kassabov and Chen Meiri for insightful
discussions, and an anonymous referee for his careful reading and helpful comments.
% Here is the bibliography.
%======================================================================================================================

\noindent Einstein Institute of Mathematics, The Hebrew University of Jerusalem, \\
Jerusalem 91904, Israel

\noindent \textit{Current address:} Department of Mathematics, Weizmann Institute of Science,\\
Rehovot 76100, Israel\\
and\\
Computer Science Division, The Open University of Israel,\\
Raanana 43107, Israel\\
\textit{E-mail address: }uzy.hadad@gmail.com


\begin{thebibliography}{12}


\bibitem[ALW]{ALW} N. Alon, A. Lubotzky and A. Wigderson, Semi-direct
    Product in Groups and Zig-Zag Product in Graphs: Connection and
    applications. 42nd IEEE Symposium on Foundation of Computer Science
    (Las Vegas, NY, 2001), 630-637, IEEE Computer Soc. Los Alamitos, CA,
    2001.


\bibitem[And]{And}S. Andreadakis, On the Automorphisms of Free Groups
    and Free Nilpotent Groups. Proc. London Math. Soc. (3) 15 1965
    239--268.

\bibitem[AR]{AR} N. Alon and Y. Roichman, Random Cayley Graphs and
    Expanders. Random Structures Algorithms 5 (1994), no. 2, 271--284.


\begin{comment}
\bibitem[Bag]{Bag} E. Bagno, Kazhdan Constants of some Colored
    Permutation Groups. J. Algebra 282 (2004), no. 1, 205--231.
    

\bibitem[BH]{BH} R. Bacher and P. de la Harpe, Exact Values of Kazhdan
    Constants for some Finite Groups.  J. Algebra 163 (1994), no. 2,
    495--515.
\end{comment}

\bibitem[BHV]{BHV} B. Bekka, P. de la Harpe and A. Valette, Kazhdan's Property ($T$).
   New Mathematical Monographs, 11. Cambridge University Press, Cambridge, 2008.

\bibitem[BNV]{BNV}
 S. R. Blackburn, P. M. Neumann, and G. Venkataraman, Enumeration of Finite Groups.
 Cambridge Tracts in Mathematics, 173. Cambridge University Press, Cambridge, 2007.

\bibitem[Bur]{Bur} M.  Burger,  Kazhdan Constants for ${\rm SL}(3,Z)$. J.
    Reine Angew. Math. 413 (1991), 36--67.


\bibitem[For]{For}E. Formanek, Fixed Points and Centers of Automorphism
    Groups of Free Nilpotent Groups. Comm. Algebra 30 (2002), no. 2,
    1033--1038.

\bibitem[HB]{HB} B. Huppert and N. Blackburn, Finite Groups II.
    Grundlehren der Mathematischen Wissenschaften [Fundamental
    Principles of Mathematical Sciences], 242. AMD, 44. Springer-Verlag,
    Berlin-New York, 1982.

\bibitem[HLW]{HLW} S. Hoory, N. Linial, and A. Wigderson, Expander
    Graphs and their Applications. Bull. AMS,
    43(2006) 439--561.

\bibitem[HV]{HV} P. de la Harpe and A. Valette, La propri\'{e}t\'{e} (T) de
    Kazhdan pour les Groups Localement Compacts. Ast\'{e}risque,
    175 (1989).

\bibitem[GZ]{GZ} T. Gelander and A. \.Zuk, Dependence of Kazhdan Constants
    on Generating Subsets. Israel J. Math. 129 (2002), 93--98.

\bibitem[Ge]{Ge}S. M. Gersten, A Presentation for the Special
    Automorphism Group of a Free Group. J. Pure Appl. Algebra 33 (1984),
    no. 3, 269--279.

%\bibitem[Hel]{Hel}G. T. Helleloid, Automorphism Groups of Finite
%    p-Groups: Structure and Applications.  PhD thesis, The University of Texas at
%    Austin, 2007.


\bibitem[Kas1]{Kas1} M. Kassabov, Universal Lattices and Unbounded Rank Expanders.
Invent. Math. 170 (2007), no. 2, 297--326.

\begin{comment}
\bibitem[Kas2]{Kas2} M. Kassabov, Kazhdan Constants for ${\rm SL}\sb
    n({\mathbb {Z}})$. Internat. J. Algebra Comput. 15 (2005), no. 5-6,
    971--995.
\end{comment}

\bibitem[Kaz]{Kaz} D.A. Kazhdan, On the Connection of the Dual Space of a
    Group with the Structure of its Closed Subgroups. Funk. Anal. Pril. 1 (1967) 71--74.

\bibitem[KLN]{KLN} M. Kassabov, A. Lubotzky and N. Nikolov, Finite Simple
    Groups as Expanders.  Proc. Natl. Acad. Sci. USA 103 (2006), no. 16,
    6116--6119.

\bibitem[LR]{LR} Z. Landau and A. Russell, Random Cayley Graphs are
    Expanders: A Simple Proof of the Alon-Roichman Theorem. Electron.
    J. Combin. 11 (2004), no. 1.

\bibitem[LS]{LS} P. Loh and L. Schulman, Improved Expansion of Random Cayley
    Graphs. Discrete Mathematics and Theoretical Computer science, 6,
    (2004), 523-528.

\bibitem[Lub]{Lub} A. Lubotzky, Discrete Groups, Expanding Graphs and
    Invariant Measures. Birkh\"{a}user, 1994.

\bibitem[LP]{LP} A. Lubotzky and I. Pak, The Product
    Replacement Algorithm and Kazhdan's Property (T). J. Amer. Math. Soc.
    14 (2001), no. 2, 347--363.

\bibitem[LW]{LW} A. Lubotzky and B. Weiss, Groups and Expanders. Expanding
    graphs (Princeton, NJ, 1992), 95--109, DIMACS Ser. Discrete Math.
    Theoret. Comput. Sci., 10, Amer. Math. Soc., Providence, RI, 1993.


\bibitem[MAK]{MAK} K. Mallahi-Karai, Relative Growth and Kazhdan Property
    for Arithmetic Groups. PhD thesis Yale univeristy 2006.


\bibitem[MW]{MW} R. Meshulam and A. Wigderson, Expanders in Group Algebras.
    Combinatorica 24 (2004), no. 4, 659--680.

\bibitem[Ne]{Ne} H. Neumann, Varieties of Groups. Springer-Verlag
    New York, Inc., New York, 1967.
\begin{comment}
\bibitem[PZ]{PZ} I. Pak and A. \.Zuk, On Kazhdan Constants and
    Mixing of Random Walks. Int. Math. Res. Not. (2002), no. 36, 1891--1905.
\end{comment}

\bibitem[RVW]{RVW} O. Reingold, S. Vadhan and A. Wigderson, Entropy Waves, the Zig-Zag Graph Product, and New Constant-Degree Expanders,
Annals of Mathematics, vol. 155, no.1, pp. 157-187, 2002 (preliminary version appeared in the proceedings

\bibitem[Rot]{Rot}  J. J. Rotman,  An Introduction to the
    Theory of Groups. Fourth edition. Graduate Texts in Mathematics, 148.
    Springer-Verlag, New York, 1995.

\bibitem[RSW]{RSW} E. Rozenman, A. Shalev and A. Wigderson,
Iterative construction of Cayley expander graphs
Theory of Computing, vol. 2, no. 5, pp. 91-120, 2006 (preliminary version appeared in the proceedings of STOC 2004).

\bibitem[Sh1]{Sh1} Y. Shalom, Bounded Generation and Kazhdan Property (T).
    Publ. Math. IHES, 90 (1999), 145--168.


\end{thebibliography}
 \end{document}